\documentclass[12pt]{amsart}
\usepackage{amssymb}
\usepackage{latexsym}
\usepackage{amsfonts}
\newtheorem{theorem}{Theorem}[section]
\newtheorem{lemma}[theorem]{Lemma}
\newtheorem{corollary}[theorem]{Corollary}
\newtheorem{proposition}[theorem]{Proposition}

\newtheorem{mainthm}{Theorem}
\newtheorem{maincor}{Corollary}

\newtheorem{definition}[theorem]{Definition}

\newtheorem{remark}[theorem]{Remark}

\newcommand{\ble}{\begin{lemma}}
\newcommand{\ele}{\end{lemma}}
\newcommand{\bth}{\begin{theorem}}
\renewcommand{\eth}{\end{theorem}}
\newcommand{\bpr}{\begin{proposition}}
\newcommand{\epr}{\end{proposition}}
\newcommand{\bco}{\begin{corollary}}
\newcommand{\eco}{\end{corollary}}
\newcommand{\bde}{\begin{definition}}
\newcommand{\ede}{\end{definition}}
\newcommand{\beq}{\begin{equation}}
\newcommand{\eeq}{\end{equation}}
\newcommand{\bpf}{\begin{proof}}
\newcommand{\epf}{\end{proof}}

\newcommand{\into}{\hookrightarrow}
\newcommand{\sbe}{\subseteq}
\newcommand{\De}{\Delta}
\newcommand{\bA}{{\bf A}}
\newcommand{\bF}{{\bf F}}
\newcommand{\cC}{{\mathcal C}}
\newcommand{\cE}{{\mathcal E}}
\newcommand{\cK}{{\mathcal K}}
\newcommand{\Aut}{\mathop{\rm Aut}\nolimits}
\newcommand{\ad}{\mathop{\rm ad}\nolimits}
\newcommand{\diag}{\mathop{\rm diag}\nolimits}
\newcommand{\id}{\mathop{\rm id}\nolimits}
\newcommand{\orG}{\overrightarrow{\Gamma}}
\newcommand{\orE}{\overrightarrow{E}}
\newcommand{\SL}{\mathop{{\rm SL}}\nolimits}
\newcommand{\PSL}{\mathop{{\rm PSL}}\nolimits}
\newcommand{\ZZ}{{\mathbb Z}}
\newcommand\Gm{\Gamma}
\newcommand{\fk}{{\rm k}}
\newcommand{\dfn}{\em}

\newcommand{\after}{\mathbin{ \circ }}
\newcommand{\GL}{\mathop{\rm GL}\nolimits}
\newcommand{\Sp}{\mathop{\rm Sp}\nolimits}
\newcommand{\SU}{\mathop{\rm SU}\nolimits}
\newcommand{\GamL}{\mathop{\rm \Gamma L}}
\newcommand{\PGaL}{\mathop{\rm P\Gamma L}}
\newcommand{\AI}[1]{\item[\rm{(#1)}]}
\newcommand{\sE}{{\mathsf{E}}}
\newcommand{\sF}{{\mathsf{F}}}
\newcommand{\tA}{{\widetilde{A}}}
\newcommand{\vep}{\varepsilon}

\newcommand\cA{{\mathcal A}}

\newcommand{\im}{{\mathop{\rm Im}}}

\newenvironment{frontmatter}{\pagestyle{empty}}{\pagestyle{plain}}
\newenvironment{keyword}{\noindent Keywords:}{}
\newcommand{\ead}{\email}

\newcommand{\corref}[1]{}
\newcommand{\sep}{\ }
\begin{document}
\begin{frontmatter}
\title{A classification of Curtis-Tits amalgams}
\author{Rieuwert J. Blok \corref{cor}}
\ead{blokr@member.ams.org}
\address{Department of Mathematics and Statistics,
Bowling Green State University,
Bowling Green, OH 43403,
U.S.A.}

\author{Corneliu G. Hoffman}
\ead{C.G.Hoffman@bham.ac.uk}
%\cortext[cor]{Corresponding author.}
\address{University of Birmingham,
Edgbaston, B15 2TT,
U.K.}

\thanks{The authors wish to thank the Banff International Research Station, Banff, Canada for a very productive visit as part of their Research In Teams program (08rit130).}
\begin{abstract}
A celebrated theorem of Curtis and Tits on groups with finite BN-pair shows that roughly speaking these groups are determined by their local structure.
This result was later extended to Kac-Moody groups by P.~Abramenko and  B.~M\"uhlherr. Their theorem states that a Kac-Moody group $G$ is the universal completion of an amalgam of rank two (Levi) subgroups, as they are arranged inside $G$ itself.

%Taking this result as a starting point, we define Curtis-Tits structures purely group theoretically as amalgams of groups that resemble, but are more general than, those in the Curtis-Tits theorem. There is no a priori reference to an ambient group. Indeed, there is no a priori guarantee that the amalgam will not collapse.
Taking this result as a starting point, we define a Curtis-Tits structure over a given diagram to be an amalgam of groups 
such that the sub-amalgam corresponding to a two-vertex sub-diagram is the Curtis-Tits amalgam of some rank-$2$ group of Lie type.
There is no a priori reference to an ambient group, nor to the existence of an associated (twin-) building. Indeed, there is no a priori guarantee that the amalgam will not collapse.

We then classify these amalgams up to isomorphism.
In the present paper we consider triangle-free simply-laced diagrams. 
Instead of using Goldschmidt's lemma, we introduce a new approach by applying Bass and Serre's theory of graphs of groups.
The classification reveals a natural division into two main types: "orientable" and "non-orientable" Curtis-Tits structures.
Our classification of orientable Curtis-Tits structures naturally fits with the classification of all locally split Kac-Moody groups using Moufang foundations.
In particular, our classification yields a simple criterion for recognizing when Curtis-Tits structures give rise to Kac-Moody groups.
The class of non-orientable Curtis-Tits structures is in some sense much larger. Many of these amalgams turn out to have non-trivial interesting completions inviting further study.
\end{abstract}
\maketitle
\end{frontmatter}
%\maketitle
\section{Introduction}
Kac-Moody Lie algebras are infinite dimensional Lie algebras  defined by relations analogous to the Serre relations for  finite dimensional semi-simple Lie algebras. They have been introduced in the mid sixties by V. Kac and R. Moody.  The affine Kac-Moody and generalized Kac-Moody Lie algebras have  extensive applications to theoretical  physics, especially conformal field theory, monstrous moonshine and more.

  Finite dimensional semi-simple Lie algebras admit Chevalley bases which allow the construction of Chevalley groups,  Lie-type  groups over arbitrary fields.  By analogy, J. Tits  defined Kac-Moody groups to be groups with a twin-root datum, which implies that  they are symmetry groups of Moufang twin-buildings (see \cite{Tit87,Ti1992}). In the case that the corresponding diagram is spherical, the corresponding group is a Chevalley group. These and other similar groups  play a very important role in various aspects of geometric group theory. In particular, they provide examples of infinite simple groups (see for example \cite{BurMoz00,BurMoz00a,CapRem09,Hee91}).

  Finite dimensional semisimple Lie algebras admit Chevalley bases which allow the construction of Chevalley groups,  Lie-type  groups over arbitrary fields.  By analogy, J. Tits  defined Kac-Moody groups to be groups with a twin-root datum, which implies that  they are symmetry groups of Moufang twin-buildings (see \cite{Tit87} \cite{Ti1992}). In the case that the corresponding diagram is spherical, the corresponding group is a Chevalley group. These and other similar groups  play a very important role in various aspects of geometric group theory. In particular, they provide important examples of infinite simple groups (see for example \cite{BurMoz00,BurMoz00a,CapRem09,Hee91}).

 A celebrated theorem of Curtis and Tits on groups with finite BN-pair  shows that roughly speaking these groups are determined by their local structure, that is by an amalgam of rank two algebraic groups. This theorem was later extended by Timmesfeld (see \cite{Tim98,Tim03,Tim04,Tim06} for spherical groups and by P. Abramenko and  B.~M\"uhlherr in \cite{AbrMuh97} to 2-spherical Kac-Moody groups.

This theorem states that the Kac-Moody groups are the universal completion of the concrete amalgam of their Levi subgroups. 
 In case that the amalgam is unique, this suffices to recognize the group. In general however, this is an inconvenience since it is usually easy to recognize isomorphism classes of subgroups put perhaps not so easy to globally manage their embedding. This is the reason that one often restricts to the so called  "split" Kac-Moody groups, that is, groups in which the embedding is the natural one. However "twisted" versions of Kac-Moody groups do exist, as constructed in ~\cite{Hee91,Ram95} and they in turn give Curtis-Tits amalgams.

%In the paper~\cite{Ti1992},  J.~Tits considered groups of Kac-Moody type and their twin-buildings and conjectured that, under certain conditions, these can be classified by Moufang foundations.
%He then proceeded to classify all Moufang foundations with simply laced Dynkin diagram without triangles,
%and stated that, under certain local conditions  these are in bijection with the maps from the fundamental group of the diagram to the automorphism group of the field in question. The correspondence between foundations and Kac-Moody groups  was later proved by B.~M\"uhlherr in~\cite{Mu1999}.  We realized that  while any Moufang foundation gives rise to a Curtis-Tits amalgam, the converse is not necessarily true.

A natural question is therefore the following: how can one recognize these amalgams as abstract group amalgams? More generally  one would like to  classify all amalgams that are "locally" isomorphic to the usual Curtis-Tits ones and identify their  universal completions.
 In this paper we use a variation of Bass-Serre theory to classify all Curtis-Tits structures over simply laced diagrams without triangles. As a by-product we obtain a description of all Kac-Moody groups in this case.

\medskip

 Throughout the paper $\fk$ will be a commutative  field of order at least $4$. We need the restriction on the order for the classification of the amalgams.
Precise definitions will be given in Section~\ref{section:ct-structures}.
A Curtis-Tits (CT) structure over $\fk$ with (simply laced) Dynkin diagram $\Gamma$ over a finite set $I$ is an amalgam $\cA=\{G_i,G_{i,j}\mid i,j\in I\}$ whose rank-$1$ groups $G_i$ are isomorphic to $\SL_2(\fk)$, where
 $G_{i,j}=\langle G_i,G_j\rangle$, and in which $G_i$ and $G_j$ commute if $\{i,j\}$ is a non-edge in $\Gamma$ and
   are embedded naturally in $G_{i,j}\cong \SL_3(\fk)$ if $\{i,j\}$ is an edge in $\Gamma$.
We are only interested in CT structures that admit a non-trivial completion. The universal completion of a (non-collapsing) Curtis-Tits structure is called a Curtis-Tits group.

%It is demonstrated that the amalgam uniquely determines a subgroup $T$ that is locally a torus.
%We call the Curtis-Tits structure orientable if it is possible to choose a family $\{X_i\}_{i\in I}$ of root groups,
% $X_i\le G_i$ such that $X_i$ and $X_j$ are in the same Borel subgroup of $G_{i,j}$ with respect to $G_{i,j}\cap T$, for every $i,j\in I$.
%Examples of ``orientable'' Curtis-Tits (OCT) structures are obtained by applying the Curtis-Tits theorem
% to locally split Kac-Moody groups $G$.
%The amalgam gives a presentation of $G$ as the completion of the  rank-$1$ and rank-$2$ Levi components of standard parabolic subgroups relative to a twin BN-pair of type $\Gamma$ for $G$.
%
%We prove that orientable Curtis-Tits structures are exactly the ones coming from Moufang foundations. However, there exist non-orientable Curtis-Tits structures.

In fact, a slight extension of our methods allows to classify Curtis-Tits structures for a larger class of diagrams, including for instance all 3-spherical Dynkin diagrams. However, in order to present these new methods and new results in a transparent manner, we chose to restrict to all simply-laced diagrams without triangles, just as Tits did in his classification of Moufang foundations for these diagrams in~\cite{Ti1992}.

Curtis-Tits groups that are not Kac-Moody groups do exist. As an application in~\cite{BloHof09b} we  give constructions of all possible Curtis-Tits structures with diagram $\tA_n$, realizing them as concrete amalgams inside their respective non-trivial completions. This leads us to describe two very interesting collections of groups. The first is a collection of twisted versions of the Kac-Moody group $\SL_n(\fk[t,t^{-1}])$ whose natural quotients are labeled by the cyclic algebras of center $\fk$. The corresponding twin-building is related to Drinfeld's vector bundles over a non-commutative projective line. The second is  a collection of Curtis-Tits groups that are not Kac-Moody groups. One of these maps surjectively to $\Sp_{2n}(q)$, $\Omega_{2n}^{+}(q)$, and $\SU_{2n}(q^l)$, for all $l\ge 1$, making this family of unitary groups into a family of expanders. 
%This Curtis-Tits group also admits a canonical quantum Clifford algebra that seems to have close relations to q-CCR algebras considered by physicists. 

In the case of hyperbolic diagrams we hope to be able to prove that all the resulting groups are simple. Of course those that are Kac-Moody groups are  simple by the results of Caprace and Remy \cite{CapRem09}.

\medskip
 Our main result is the following.

 \begin{mainthm}\label{mainthm:CT structures} Let $\Gamma$ be a simply laced Dynkin diagram with no triangles and $\fk$ a field with at least 4 elements.
There is a natural bijection between isomorphism classes of CT-structures over the field $k$ on a graph $\Gamma$ and elements of the set $\{\Phi\colon\pi(\Gamma, i_0)\to \ZZ_2\times\Aut(\fk) |\ \Phi \mbox{ is a group homomorphism}\}$
 \end{mainthm}

Here, $\pi(\Gamma,i_0)$ denotes the fundamental group of the graph $\Gamma$ with base point $i_0$.
As mentioned above, the  motivation for the work came from the Curtis-Tits amalgam presentations for Kac-Moody groups. In fact in the spherical case these were proved to be the only such amalgams. Surprisingly, in general they form a small minority of all amalgams. More precisely they are those amalgams in the theorem corresponding to maps $\Phi$ so that $\im(\Phi) \le \Aut(\fk)$. We call such amalgams "orientable". The relation between Kac-Moody groups and orientable CT amalgams is made via Moufang foundations. By results of Tits~\cite{Ti1992} and M\"uhlherr~\cite{Mu1999}, Moufang foundations of type
 $\Gamma$ over $\fk$ are classified by homomorphisms  from $\pi(\Gamma,i_0)$ to $\Aut(\fk)$.
Moreover, by the main result of M\"uhlherr~\cite{Mu1999}, any
foundation with a simply laced diagram in which every $A_2$-residue
is of type $A_2(\fk)$ (i.e.~locally split) can be ``integrated''. Combining these results with Theorem~\ref{mainthm:CT structures}, we can then prove the following corollary:

\begin{maincor}\label{mainthm:OCT=Moufang foundation}
Let $\Gamma$ be a simply laced Dynkin diagram with no triangles and $\fk$ a field with at least 4 elements.
The universal completion of a  Curtis-Tits structure over a commutative field $\fk$ and diagram $\Gamma$ is a locally split Kac-Moody group over $\fk$ with Dynkin diagram $\Gamma$ (and $\cA$ is the Curtis-Tits amalgam for this group) if and only if $\cA$ is orientable.
\end{maincor}

Note that for example in~\cite{AbrMuh97,Cap07,Ti1992} 
 the amalgam is required to live in the corresponding Kac-Moody
group.  This is rather inconvenient since it gives no intrinsic
description of the amalgam. Our result above defines Kac-Moody
groups as universal completions of certain abstract amalgams hence
giving concrete presentations for those groups. In particular,  we
can refine  Corollary~\ref{mainthm:OCT=Moufang foundation} as follows. See
Section~\ref{sec:Sound Moufang foundations and orientable
CT-amalgams} for the exact definitions.
\begin{maincor}\label{mainthm:presentation}
Let $\Gamma$ be a simply laced Dynkin diagram with no triangles and $\fk$ a field with at least 4 elements. Any locally split Kac-Moody group over $\fk$ with diagram $\Gamma$ can be defined by a twist $(\orG, \delta)$ of the corresponding split Kac-Moody group. Moreover any two twists are equivalent if they have the same fundamental group.
\end{maincor}

The corollaries above could be proved directly from the above
mentioned results of Tits and M\"uhlherr. To our knowledge however
there is no explicit correspondence in the literature to this
effect. Moreover, in the absence of Theorem~\ref{mainthm:CT
structures}, it is not immediately obvious that different choices of
an orientable CT amalgam would give different foundations.  See also
Corollary~\ref{cor:classification of KM groups} for a more precise
construction of the amalgams in the spirit of \cite{Cap07} (see the
application to Theorem A in loc.~cit.).

The paper is organized as follows. In Section~\ref{section:ct-structures} we define Curtis-Tits structures, morphisms and prove some general technical lemmas. In Section~\ref{sec: Bass serre} we introduce our modification of Bass-Serre theory and prove Theorem~\ref{mainthm:CT structures}. In Section~\ref{section:Curtis-Tits Theorem} we prove Corollary~\ref{mainthm:OCT=Moufang foundation}
 and in Section~\ref{section:twists} we prove Corollary~\ref{mainthm:presentation}.

\medskip
\paragraph{\bf Acknowledgement}
%This paper was started and planned while the authors enjoyed the hospitality of the Banff International Research Station on a Research in Teams program. We are extremely grateful to  BIRS for the pleasant and highly productive stay.
The final version of the paper was written during some wonderful, if very claustrophobic and accident prone three weeks in Birmingham. We thank Irina and Karin for putting up with it all.

\section{CT-structures}\label{section:ct-structures}
In this section we introduce the notion of a CT-structure over a commutative field and define its category.
Throughout the paper $\fk$ will be a commutative field.

\bde\label{dfn:standard pair}
Let $V$ be a vector space of dimension $3$ over $\fk$.
We call $(S_1,S_2)$ a {\em standard pair} for $S=\SL(V)$ if
 there are decompositions $V=U_i\oplus V_i$, $i=1,2$, with $\dim(U_i)=1$
  and $\dim(V_i)=2$ such that $U_1\sbe V_2$ and $U_2\sbe V_1$ and
   $S_i$ centralizes $U_i$ and preserves $V_i$.

One also calls $S_1$ a {\dfn standard complement} of $S_2$ and vice-versa.
We set $D_1= N_{S_1}(S_2)$ and $D_2=N_{S_2}(S_1)$.
A simple calculation shows that $D_i$ is a maximal torus in $S_i$, for $i=1,2$.
In general if $G\cong \SL_3(\fk)$, then $(G_1,G_2)$ is a standard pair for $G$ if there is an isomorphism
 $\psi\colon G\to S$ such that $\psi(G_i)=S_i$ for $i=1,2$.
\ede

\bde\label{dfn:standard basis for standard pair}
Given a standard pair $(S_1,S_2)$, a {\em standard basis} for $(S_1, S_2)$ is an ordered basis $\sE_0=(e_1,e_2,e_3)$ of $V$ such that $V_1=\langle e_1,e_2\rangle$, $U_1= \langle e_3\rangle $,  $U_2=\langle e_1\rangle$,  and  $V_2=\langle e_2,e_3\rangle$.
\ede
Identifying $S$ with $\SL_3(\fk)$ via its left action on $V$ with respect to $\sE_0$, yields
$$\begin{array}{@{}lll}
 S_1=\left\{\left.\left(\begin{array}{@{}cc@{}}A & 0 \\0 & 1\end{array}\right) \right| A \in \SL_{2}(\fk)\right\}
& \mbox{  and } &
S_2=\left\{\left.\left(\begin{array}{@{}cc@{}}1 & 0 \\0 &  A\end{array}\right) \right|A \in \SL_2(\fk)\right\}\\
\mbox{so that }  &&\\
 D_1 =\left\{\left.\left(\begin{array}{@{}ccc@{}}
a & 0 & 0 \\
0 & a^{-1} & 0\\
0 & 0 & 1 \end{array}\right) \right|a \in \fk^*\right\}
& \mbox{  and } &
D_2=\left\{\left.\left(\begin{array}{@{}ccc@{}}
1 & 0 & 0\\
0 & a  & 0\\
0 & 0 & a^{-1} \end{array}\right) \right|a \in \fk^*\right\}
\end{array}$$
\ble\label{lem:two complements} Let $S_1$ and $S_2$ be a standard
pair for $S=\SL_3(\fk)$, where $\fk$ has at least four elements.
Then $S_1$ has exactly one standard complement $S_2'\ne S_2$
normalized by $D_1$. \ele \bpf Since $\fk$ has at least four
elements, $D_1$ uniquely determines three $1$-dimensional
eigenspaces
 and $S_1$ fixes all vectors in exactly one of these eigenspaces, say $E_1$.
In the notation above, these are $E_1=U_1$, $U_2$ and $V_1\cap V_2$.
Thus any standard complement $S_2$ to $S_1$ that is normalized by
$D_1$ is completely determined by the eigenspace $E \ne E_1$ that it
fixes vector-wise. There are two choices. \qed\epf

 We will need the following lemma.
 \ble\label{lem:standard pair of tori}
With the notations above, $D_{1}=C_{S_{1}}(D_2)$ and  $D_{2}=C_{S_{2}}(D_1)$. Moreover, $D_{2}$ is the only torus in $S_{2}$ that is normalized by $D_{1}$.
 \ele
 \bpf Note that if $T$ is a torus in $S_{2}$ then  $N_{S}(T)$ is the set of monomial matrices so $N_{S_{1}}(T)$ only contains one torus which is $C_{S_{1}}(T)$. The conclusion follows.  \qed\epf

\bde\label{dfn:diagram}
A {\em simply laced Dynkin diagram} over the set $I$ is a simple graph $\Gamma=(I,E)$. That is, $\Gamma$ has vertex set $I$, and an edge set $E$ that contains no loops or double edges.
\ede

\bde\label{dfn:amalgam}
An {\em amalgam} over a set $I$ is a collection $\cA=\{G_i,G_{i,j}\mid i,j\in I\}$ of groups, together with a collection $\varphi=\{\varphi_{i,j}\mid i,j\in I\}$ of monomorphisms $\varphi_{i,j}\colon G_i\into G_{i,j}$, called {\em inclusion maps}.
A {\em completion} of $\cA$ is a group $G$ together with a collection  $\phi=\{\phi_i,\phi_{i,j}\mid i,j\in I\}$ of homomorphisms $\phi_i\colon G_i\to G$ and $\phi_{i,j}\colon G_{i,j}\to G$, such that for any $i,j$ we have $\phi_{i,j}\after\varphi_{i,j}=\phi_i$.
For simplicity we denote by $ \bar G_i=\varphi_{i,j}(G_i)\le G_{i,j}$.
The amalgam $\cA$ is {\em non-collapsing} if it has a non-trivial completion.
A completion $(\hat{G},\hat{\phi})$ is called {\em universal} if for any completion $(G,\phi)$ there is a unique surjective group homomorphism $\pi\colon \hat{G}\to G$ such that $\phi=\pi\after\hat{\phi}$.
\ede

\bde\label{dfn:weak CT structure}
Let $\Gm=(I,E)$ be a simply laced Dynkin diagram.
A {\dfn Curtis-Tits structure over $\Gm$} is a non-collapsing amalgam $\cA(\Gm)=(G_{i},G_{i\, j} | i, j \in I)$ such that
\begin{enumerate}
  \item[(CT1)]  for any vertex $i$, the group $ G_i = \SL_2(\fk)$ and for each pair $i,j \in I$,
  $$G_{i,j}\cong\left\{\begin{array}{cc}\SL(V_{i,j}) & \mbox{if} \ \{i,j\}\ \in E \\ G_i\ast G_j & \mbox{if} \ \{i,j\}\ \not\in E\end{array}\right.,$$
  where $V_{i,j}$ is a 3-dimensional vector space over $\fk$
   and $\ast$ denotes central product;
  \item[(CT2)] if $\{i,j\}\in E$ then $(\bar G_{i}, \bar G_{j})$ is a standard pair in $G_{i,j}$.
\end{enumerate}
\ede

\bde
A Dynkin diagram is {\dfn admissible} if it is connected and has no circuits of length $\le 3$.
\ede

From now on $\Gamma=(I,E)$ will be an admissible Dynkin diagram and $\cA=\cA(\Gamma)=\{G_i,G_{i,j}\mid i,j\in I\}$ will be a non-collapsing Curtis-Tits structure over $\Gm$.

 \begin{lemma} \label{lem:invariant di}
If the Dynkin Diagram is admissible and  $i, j, k$ are vertices such that $\{i,j\}$and $\{j,k\}$ are edges then $N_{G_{i,j}}(\bar G_{i})\cap \bar G_{j}= N_{G_{jk}}(\bar G_{k})\cap \bar G_{j}$\end{lemma}

\bpf (See also \cite{D05}) Let $(G,\phi)$ be a non-trivial
completion of $\cA$ and identify $\cA$ with its $\phi$-image in $G$.
Let $\bar D_{j}^{i}=N_{G_{i,j}}(\bar G_{j})\cap \bar G_{i}$. It
follows from the fact that the  nodes $i$ and $k$ are not connected
in $\Gamma$ that $\bar D^{i}_{j}$ commutes with $\bar D^{k}_{j}$.
Note that if $g\in \bar D^{i}_{j}$ then $(\bar D^{j}_{k})^{g}$
commutes with $\bar D^{k}_{j}$ so  $\bar D^{i}_{j}$ is a torus that
normalizes the torus $\bar D^{j}_{k}$ of $\bar{G}_j$. By
Lemma~\ref{lem:standard pair of tori}, $\bar D^i_j$ only normalizes
$\bar D^j_i$  and so $\bar  D^{j}_{k}=\bar  D^{j}_{i}$. \qed\epf

Lemma~\ref{lem:invariant di} motivates the following definition.

\bde\label{dfn:D_i}
For $i,j\in I$, we let $\bar D_i=N_{G_{i,j}}(\bar G_{j}) \cap \bar  G_{i}$, where $\{i,j\}\in E$.
Note that this defines $\bar  D_{i}$ for all $i$ since $\Gamma$ is connected. We also denote $D_i =\varphi _{i,j}^{-1}(\bar  D_i)$.
\ede

As we saw after Definition~\ref{dfn:standard pair}, for each $i\in
I$, the group $\bar  D_{i}$ is a torus in $ \bar G_{i}$.
Lemma~\ref{lem:invariant di} allows us to glue tori together.
\ble\label{lem:ij-torus} If $\{i,j\}\in E$, then $\bar  D_i$ and $
\bar D_j$ are contained in a unique common maximal torus $D_{i,j}$
of $G_{i,j}$. \ele \bpf Clearly in any completion of the amalgam,
both $ \bar D_i$ and $ \bar D_j$ normalize  $\bar G_i$ and $\bar
G_j$ so we have $\bar D_i, \bar  D_j\le N_{G_{i,j}}(\bar G_i)\cap
N_{G_{i,j}}(\bar G_j)=D_{i,j}$, which is the required maximal torus.
\qed\epf

\bde\label{dfn:OCT}
Note that a torus in $\SL_{2}(\fk)$ uniquely determines a pair of opposite root groups $X_{+}$. and $X_{-}$. We now choose one root group  $X_{i}$ normalized by the torus $D_{i}$ of $G_{i}$ for each $i$.
An {\it orientable Curtis-Tits (OCT) structure} (respectively orientable Curtis-Tits (OCT) group) is a CT structure that admits a system $\{X_i\mid i\in I\}$ of root groups as above such that for any $i,j\in I$, the groups $\varphi _{i,j} (X_{i})$ and
$\varphi_{j,i}(X_{j})$ are contained in a common Borel subgroup $B_{i,j}$ of $G_{i,j}$.
\ede

\subsection{Morphisms}
In this subsection, for $k=1,2$, let $\Gm^k=(I^k,E^k)$ be a Dynkin diagram.

\bde\label{dfn:Dynkin hom}
 A {\dfn homomorphism} between the Dynkin diagrams $\Gm^1$ and $\Gm^2$ is a map $\gamma\colon I^1\to I^2$
such that for any $i,j\in I$ with $\{i,j\}\in E^1$ also $\{\gamma(i),\gamma(j)\}\in E^2$.
We call $\gamma$ an isomorphism if $\gamma$ is bijective and $\gamma^{-1}$ is also a homomorphism of Dynkin diagrams, that is $\{i,j\}\in E^1$ if and only if
$\{\gamma(i),\gamma(j)\}\in E^2$ for all $i,j\in I^1$.
We call $\gamma$ an automorphism if $\gamma$ is an isomorphism and $\Gm^1=\Gm^2$.
\ede

\noindent
Now, for $k=1,2$, let $\cA^k=\{G_i^k,G_{i,j}^k\mid i,j\in I^k\}$ be a CT structure with admissible Dynkin diagram $\Gm^k$.

\bde\label{dfn:CT hom}
A {\dfn homomorphism} between the amalgams $\cA(\Gm^1)$ and $\cA(\Gm^2)$ is a
pair $(\gamma,\phi)$, where $\gamma\colon \Gm^1\to\Gm^2$ is a homomorphism and
$\phi=\{ \phi_i, \phi_{i,j}\mid i,j\in I^1\}$ where $\phi\colon G_i^1\to  G_{\gamma(i)}^2$ and $\phi_{i,j}  \colon G_{i,j}^1\to G_{\gamma(i),\gamma(j)}^2$are group homomorphisms  such that $$\phi_{i,j}\after \varphi_{i,j}^1 =\varphi_{{\gamma(i)},{\gamma(j)}}^2\after  \phi_i.$$
We call $(\gamma,\phi)$ an {\dfn isomorphism} of amalgams if $\gamma$ is an isomorphism, $\phi_i$ and $\phi_{i,j}$ are bijective for all $i,j\in I$, and $(\gamma^{-1},\phi^{-1})$ is a homomorphism of amalgams. Note that, if $(\gamma,\phi)$ is an isomorphism, we can relabel the elements of $\Gamma^2$ and assume $\gamma =\id$. For most of the following we will do so and denote the isomorphism simply by $\phi$.
 \ede

\ble\label{lem:cA hom on D_i}
With the notation of Definition~\ref{dfn:CT hom}, suppose that $\phi_i\colon  G_i^1\to  G_{\gamma(i)}^2$ is surjective for all $i\in I^1$. Then, the homomorphism $\phi_i$ restricts to a group homomorphism
 $\phi_i\colon D_i^1\to D_{\gamma(i)}^2$ for all $i\in I^1$.
\ele \bpf Consider any edge $\{i,j\}\in E^1$ and write
$\phi=\phi_{i,j}$ for short. Since  $\phi$ is a homomorphism with
$\phi(\bar{G_i}^1)=\bar{G}_{\gamma(i)}^2$ and
$\phi(\bar{G_j}^1)=\bar{G}_{\gamma(j)}^2$, we have
 $$\begin{array}{rll}
  \phi(\bar D_i^1)=\phi(N_{G_{i,j}^1}(\bar G_j^1)\cap \bar G_i^1)
 &\le N_{\phi(G_{i,j}^1)}(\phi(\bar G_j^1))\cap\phi(\bar G_i^1)&=N_{\phi(G_{i,j}^1)}(\bar G_{\gamma(j)}^2)\cap \bar G_{\gamma(i)}^2\\
 &\le N_{G_{\gamma(i) \gamma(j)}^2}(\bar G_{\gamma(j)}^2)\cap \bar G_{\gamma(i)}^2&=\bar D_{\gamma(i)}^2.\\
 \end{array}$$
\qed\epf \noindent

\subsection{Automorphisms of $\cA(A_2)$}
Let $W$ be a (left) vector space of dimension $n$ over $\fk$.
Let $G=\SL(W)$ act on $W$ as the matrix group $\SL_n(\fk)$ with respect to some fixed basis $\sE=\{e_i\mid i=1,2,\ldots,n\}$.
Let $\omega\in \Aut(\SL_n(\fk))$ be the automorphism
 given by $$A\mapsto {}^tA^{-1}$$
where ${}^tA$ denotes the transpose of $A$.

Let $\Phi=\{(i,j)\mid 1\le i\ne j\le n\}$.
For any $(i,j)\in \Phi$ and $\lambda\in k$, we define the {\dfn root group}
 $X_{i,j}=\{X_{i,j}(\lambda)\mid \lambda\in k\}$, where $X_{i,j}(\lambda)$ acts as
 $$\begin{array}{@{}ll} e_j\mapsto e_j+\lambda e_i & \mbox{ and }\\
                                      e_k\mapsto e_k & \mbox{ for all }k\ne j.\\
                                                   \end{array}$$
Let $\Phi_+=\{(i,j)\in \Phi\mid i<j\}$ and $\Phi_-=\{(i,j)\in \Phi\mid j<i\}$.
We call $X_{i,j}$ {\dfn positive} if $(i,j)\in \Phi_+$ and {\dfn negative} otherwise.
Let $H$ be the torus of diagonal matrices in $\SL_n(\fk)$ and for $\vep\in \{+,-\}$, let $X_\vep=\langle X_{i,j}\mid (i,j)\in \Phi_\vep\rangle$ and $B_\vep= H\ltimes X_\vep$.
\ble\label{lem:omega}
\
\begin{itemize}
\AI{a} If $n=2$, then $\omega$ is given by conjugation with
 $$\cE=\left(\begin{array}{@{}cc@{}} 0 & -1 \\ 1 & 0 \\\end{array}\right)\in\SL_2(\fk).$$

\AI{b} If $n\ge 3$, then $\omega$ cannot be represented by an element of $\GL_n(\fk)$.
\AI{c}
 $X_{i,j}^\omega=X_{j,i}$ for all $(i,j)\in \Phi$ and $B_\vep^\omega=B_{-\vep}$, for $\vep\in\{+,-\}$.
\end{itemize}
\ele \bpf (a) and (c) Straightforward. (b) If $n\ge 3$, then
$\omega$ does not even preserve eigenvalues, so it is certainly not
linear. \qed\epf

Let $\GamL_n(\fk)$ be the group of all semilinear automorphisms of the vector space $W$
 and let $\PGaL_n(\fk)=\GamL_n(\fk)/Z(\GamL_n(\fk))$.
Then $\GamL_n(\fk)\cong \GL_n(\fk)\rtimes\Aut(\fk)$, where we view $t\in\Aut(\fk)$
 as an element of $\GamL_n(\fk)$ by setting $((a_{i,j})_{i,j=1}^n)^t=(a_{i,j}^t)_{i,j=1}^n$.
The automorphism group of $\SL_n(\fk)$ can be expressed using $\PGaL_n(\fk)$ and $\omega$
 as follows \cite{SchVan28}.
\ble\label{lem:Aut(SL)}
$$\Aut(\SL_n(\fk))=\left\{\begin{array}{@{}ll} \PGaL_n(\fk) & \mbox{ if } n=2;\\
\PGaL_n(\fk)\rtimes \langle \omega\rangle & \mbox{ if }n\ge 3.\\\end{array}\right.$$
\ele

 \section{Bass-Serre theory on graphs of groups}\label{sec: Bass serre}
From a CT-structure $\cA$ we will construct a graph of groups in the sense of Bass-Serre (see \cite{Bas93,Bas96,Se80,BasLub01}). We review the relevant definitions.

\bde \label{def: Bass-Serre graphs}
Let $\Gamma=(I,E)$ be an admissible Dynkin diagram. Following \cite{Bas93} we define a directed graph $\orG=(I,\orE)$ where for each edge $\{i,j\}\in E$ we introduce directed edges $(i,j)$ and $(j,i)$ in $\orE$. For every $e\in \orE$ we denote the reverse edge by $\bar e$. Moreover we denote by $\delta_0(e)$ the starting node of the oriented edge $e$.
\ede

\bde  \label{def:graph of groups} A {\em graph of groups} is a pair $(\cC, \orG)$ where $\orG$ is a graph as above and $\cC$ associates to each $i\in I$ a group $A_i$ and to each directed edge $e\in \orE$ a group $A_e=A_{\bar e}$.
Moreover, for each vertex $i$ on a (directed) edge $(i,j)$ we have a monomorphism $\alpha_{i,j}\colon A_{i,j}\to A_i$.
\ede

\bde\label{dfn:inner morphism}
Given graphs of groups $(\cC^{(k)},\orG^{(k)})$ for $k=1,2$, an {\em inner morphism} is a pair $(\phi,\gamma)$, where
 $\gamma$ is a morphism of Dynkin diagrams and
 $\phi=\{\phi_i,\phi_{i,j}\mid i,j\in I, (i,j)\in\orE\}$ is a collection of group homomorphisms $\phi_i\colon A_i^{(1)}\to A^{(2)}_{\gamma(i)}$
  and
$\phi_{i,j}\colon A_{i,j}^{(1)}\to A^{(2)}_{\gamma(i)\gamma(j)}$
  so that for each $(i,j)\in \orE$ there exists an element
   $\delta_{i,j}\in A_{\gamma(i)}$ so that for all $s\in A_{i,j}$,
   $$\phi_{i}(\alpha_{i,j}(s))=\delta_{i,j}\alpha_{\gamma(i),\gamma(j)}(\phi_{i,j}(s))\delta_{i,j}^{-1}.$$
We call an inner morphism {\em central} if $\delta_{i,j}=1$ for all $(i,j)\in \orE$.
\ede

\noindent
Given a group $G$ and a collection of subgroups $G_1,\ldots,G_k$ let
 $\Aut_G(G_1,\ldots,G_k)$ be the subgroup of $\Aut(G)$ that stabilizes each $G_i$.
Given a monomorphism of groups $\phi\colon G\to H$, there is a corresponding homomorphism $\ad(\phi)\colon\Aut_H(\phi(G))\to \Aut(G)$ such that for any $a\in \Aut_H(\phi(G))$  we have
 $\ad(\phi)(a)=\phi^{-1}\after a\after \phi$.

Assume that  $\cA=\{G_i, G_{i,j} \mid i, j \in I \}$ a CT structure with Dynkin diagram $\Gamma=(I,E)$ and $\orG=(I,\orE)$ is the directed graph associated $\Gamma$ as in Definition~\ref{def: Bass-Serre graphs}. As we know from Lemma~\ref{lem:invariant di}, for each $i\in I$ the subgroups $\bar D_i$ and $D_i$ are well defined, hence the normal subgroup $T_i$ of diagonal automorphisms in $\Aut_{G_i}(D_i)$ is uniquely determined by $\cA$.
Similarly, for each $\{i,j\}\in E$, the normal subgroup $T_{i,j}$ of diagonal automorphisms in $\Aut_{G_{i,j}}(D_{i,j})$ is uniquely determined by $\cA$. Using Lemma~\ref{lem:Aut(SL)} one finds that
$\Aut_{G_i}(D_i)\cong T_i\rtimes (\langle\cE\rangle\times\Aut(\fk))$ and $\Aut_{G_{i,j}}(D_{i,j})\cong T_{i,j}\rtimes (\langle\omega\rangle\times\Aut(\fk))$. Note that the complements to $T_i$ and $T_{i,j}$ are both isomorphic to $\ZZ_2\times\Aut(\fk)$.

\ble\label{lem:diagonal automorphism extension} Given any collection
$\{\tau_i\in T_i\mid i\in I\}$, there exist unique automorphisms
$\tau_{i,j}\in T_{i,j}$ such that $\tau=\{\tau_i,\tau_{i,j}\mid
i,j\in I\}$ is an automorphism of $\cA$. \ele \bpf First we note
that $\bar{\tau}_i=\ad(\varphi_{i,j}^{-1})(\tau_i)$ is a diagonal
(linear) automorphism of $\bar G_i\le \bar G_{i,j}$.

If $\{i,j\}\not\in E$, then $\tau_{i,j}$ is simply the central product $\bar\tau_i\ast\bar\tau_j$.
Otherwise, suppose that with respect to some basis $\{e_1,e_2,e_3\}$ of eigenvectors for $\bar D_i$ and $\bar  D_j$ we have $\bar\tau_i=\diag\{a,b,1\}$ and $\bar\tau_j=\diag\{1,c,d\}$, then let
 $\tau_{i,j}=\diag\{ac,bc,bd\}$.
\qed\epf

\bde\label{def:base} Let $\cA=\{ G_i,G_{i,j}\mid i,j\in I\}$ be a CT structure with admissible Dynkin diagram $\Gamma=(I,E)$. A {\em basis} of $\cA$ is a collection $\{\sE_{i,j}, \sE_{j,i}\mid \{i,j\}\in E \}$ so that $\sE_{i,j}=\{e_1^{i,j},e_2^{i,j},e_3^{i,j}\}$ is a standard basis for $(\bar G_i, \bar G_j)$ in $V_{i,j}$ and $E_{j,i}$ is the same basis but the ordering is reversed. Note that $E_{i,j}$ is  stabilized by $\bar D_i$ and $\bar  D_j$. The {\em edge reversal map} is the element $\rho_{i,j}$ of $\GL(V_{i,j})$ that reverses the order of the basis $\sE_{i,j}$.
\ede

Let $\sE$ be a basis for $\cA$ as in Definition~\ref{def:base}. For each $i\in I$ let $V_i$ be a vector space with basis $\{f_1^i,f_2^i\}$ identifying $G_i=\SL_2(\fk)$.
Let $\psi_{i,j}\colon G_i\to \bar{G_i}\le G_{i,j}$ be the isomorphism induced by the linear map that takes the ordered basis $(f_1^i,f_2^i)$ to  $(e_1^{i,j},e_2^{i,j})$.
This defines a graph of groups $(\cC_0,\orG)$ in the following way.
We let $\cC_0=\{\bA_i,\bA_{i,j}\mid i,j\in I, (i,j)\in \orE\}$ where $\bA_i$ is the complement in $\Aut_{G_i}(D_i)$ to $T_i$ defined with respect to
 $\{f_1^i,f_2^i\}$ and $\bA_{i,j}$ is the complement in $\Aut_{G_{i,j}}(\bar{G_i},\bar{G_j})$ to $T_{i,j}$, defined by $\sE_{i,j}$ (See Lemmas~\ref{lem:omega}~and~\ref{lem:Aut(SL)}). Finally, we may define the map $\alpha_{i,j}\colon \bA_{i,j}\to \bA_i$ as given by the restriction of $\ad(\psi_{i,j})$ to $\bA_{i,j}$, as the following lemma shows.

\ble\label{lem:proto graph of groups} The graph of groups
$(\cC_0,\orG)$ constructed above is determined by $\fk$ and the
diagram $\orG$ up to central isomorphism (but not by the particular
amalgam $\cA$). \ele \bpf First note that the construction of
$(\cC_0,\orG)$  only involves the maps $\psi_{i,j}$, which in turn
depend uniquely on the basis $\sE$ for $\cA$ and the collection
$\sF=\{\sF_i=(f_1^i,f_2^i)\mid i\in I\}$ of bases chosen for the
$V_i$. We now show that any other choice of $\sE$ and $\sF$ merely
induces a central isomorphism between the resulting graphs of
groups. Let $\sE'$ and $\sF'$ be another choice of a basis for $\cA$
and
 the $V_i$'s and let $(\cC',\orG)$ be the resulting graph of groups.
For each $i\in I$, let $t_i\in\Aut(G_i)$ be induced by the linear map sending the ordered basis $\sF_i$ to $\sF_i'$ and for each $(i,j)\in \orE$, let
$t_{i,j}\in\Aut(G_{i,j})$ be induced by the linear map sending the ordered basis $\sE_{i,j}$ to $\sE_{i,j}'$.
Then the following diagram is commutative.
$$\begin{array}{r@{}ccc@{}l}
 & G_{i,j} & \stackrel{t_{i,j}}{\longrightarrow} & G_{i,j} & \\
 \psi_{i,j} & \uparrow & & \uparrow & \psi_{i,j}'\\
 & G_{i} & \stackrel{t_{i}}{\longrightarrow} & G_{i} & \\
\end{array}$$
Since the bases defining the complements $A_{i,j}$, $A_{i,j}'$, $A_i$ and $A_i'$ all correspond via the maps in this diagram,
 also these complements themselves correspond to each other via the adjoint maps.
This shows that the map $(\phi,\gamma)\colon \cC'\to \cC_0$, where  $\gamma$ is the identity map on $\orG$ and
 $\phi=\{\phi_{i,j}=\ad(t_{i,j}),\phi_{i}=\ad(t_{i})\mid i,j\in I, (i,j)\in \orE\}$ is a central isomorphism.
\qed\epf

For the remaining of the paper we will fix the groups $G_i$, and $G_{i,j}$ as well as the bases $\sE$ and $\sF$, which in turn fix  $\bA_i, \bA_{i,j}$, the maps $\psi_{i,j}$ and the graph of groups $\cC_0$. We note that in order to specify the CT structure $\cA$ we have to make a choice for the maps $\varphi_{i,j}$.

\bde\label{dfn:concrete amalgams}
A {\em concrete CT structure} is a CT structure $\cA=\{G_i, G_{i,j}\mid i, j \in I, (i,j)\in \orE\}$, where the groups $G_i$ and $G_{i,j}$ as well as the maps $\psi_{i,j}$ are fixed as above, and such that the inclusion maps $\varphi_{i,j}$ satisfy $\ad(\varphi_{i,j})(\bA_{i,j})=\bA_i$.
The graph of groups $\cC_0$, which is naturally associated with $\cA$, is called {\em the concrete graph of groups}.
 \ede

\noindent
Consider the concrete graph of groups $\cC_0$ and  for each $l, m\in I$ consider  $\gamma= i_0, i_1, \ldots, i_n$, a path from $l=i_0$ to $m=i_n$ in $\orG$. Define $\beta_{l,m}\colon\bA_l\to \bA_m$ by setting $\beta_{l,m}(a)=\alpha_{i_n, i_{n-1}}\after \alpha_{i_{n-1},i_n}^{-1}\after \cdots\after\alpha_{i_1,i_0}\after \alpha_{i_0,i_1}^{-1}(a)$, for each $a\in \bA_l$.

\ble\label{lem:walking on trivial paths}
 The map $\beta_{l,m}$ is independent of $\gamma$.
 \ele \bpf
 Quite immediate since, for each $(i,j)\in \orE$ the map $\alpha_{j,i}\after\alpha_{i,j}^{-1}\colon \bA_i\to \bA_j$ is
 the adjoint of the isomorphism given by the linear map $V_j\to V_i$ sending the ordered basis $\sF_j=(f_1^j,f_2^j)$
  to $\sF_i=(f_1^i,f_2^i)$.
\qed\epf

\ble\label{lem:canonical basis}
Let $\cA'=\{G_i, G_{i,j}, \varphi_{i,j}' \mid i, j \in I\}$  be a CT structure.
Then, given a collection $\{A_i\le \Aut_{G_i}(D_i)\mid i \in I\}$ of  complements to the groups of diagonal automorphisms $T_i$, there exists a basis  $\sE'=\{\sE_{i,j}'\mid \{i,j\}\in E \}$ and a collection $\{A_{i,j}\mid i,j\in I\}$ of complements to $T_{i,j}$ such that, for each $\{i,j\} \in E$, $A_{i,j}$ corresponds to $\sE_{i,j}'$ and $\ad(\varphi_{i,j}')(A_{i,j}) = A_i$. The collection $\cC=\{A_i,A_{i,j}\mid i,j\in I\}$ is unique and the bases $\sE_{i,j}'$ are unique up to multiplication by a scalar in $\mbox{Fix}(\Aut(\fk))$  \ele
 \bpf
The group $T_{i,j}$ acts regularly on the set of its complements, while acting on the corresponding bases. Two bases correspond to the same complement if and only if one is a scalar multiple of the other and that scalar is fixed by $\Aut(k)$. This proves the uniqueness part of the theorem.

For the existence we first pick a random base $\sE''$ and modify it
as follows. If $\{i,j\}\in E$ then $\sE''$ determines  $A'_{i,j}$, a
complement to $T_{i,j}$. Restriction  to $G_i$ and $G_j$ determines
complements $A_i'$ and $A_j'$ to $T_i$ and $T_j$. These are
conjugates of $A_i$ and $A_j$ under diagonal automorphisms
$\tau_i\in T_i$ and $\tau_j\in T_j$. As in the proof of
Lemma~\ref{lem:diagonal automorphism extension} there exists a
diagonal automorphism $\tau_{i,j}\in T_{i,j}$ that restricts to
$\tau_i$ and $\tau_j$. Conjugating by $\tau_{i,j}$ sends $A_{i,j}'$
to a complement $A_{i,j}$ satisfying the statement of the lemma for
the edge $\{i,j\}$, while the underlying linear map transforms the
basis $\sE''_{i,j}$ to the desired basis $\sE_{i,j}'$. \qed\epf

\bco\label{co:all amalgams are concrete} Any CT-structure
$\cA'=\{G_i, G_{i,j}, \varphi_{i,j}' \mid i, j \in I\}$ is
isomorphic to a concrete one. Moreover, the isomorphism can be taken
to be diagonal. \eco \bpf By Lemma~\ref{lem:canonical basis} there
exists a basis $\sE'$ so that the corresponding collection
$\cC=\{A_i,A_{i,j}\mid i,j\in I\}$ satisfies $A_i=\bA_i$ for all
$i\in I$, and $\ad(\varphi_{i,j}')(A_{i,j}) = \bA_i$ for all
$(i,j)\in \orE$.

We now define an isomorphism $\phi=\{\phi_i,\phi_{i,j}\mid i,j\in
I\}$ from a concrete amalgam $\cA$ to $\cA'$. Recall that $\sE$ is
the basis corresponding to the complements $\bA_{i,j}$. Now let
$\phi_{i,j}\colon G_{i,j}\to G_{i,j}$ be the isomorphism induced by
the (diagonal) linear map sending $\sE_{i,j}$ to $\sE'_{i,j}$ and
let $\phi_i=\id_{G_i}$ for all $i\in I$. Now define $\cA$ by setting
$\varphi_{i,j}=\phi_{i,j}^{-1}\after \varphi_{i,j}'\after\phi_i$.
Note that $\cA =\{G_i, G_{i,j}, \varphi_{i,j} \mid i, j \in I ,
(i,j)\in \orE\}$ is concrete since
$\ad(\varphi_{i,j})(\bA_{i,j})=\ad(\varphi_{i,j}')\after
\ad({\phi_{i,j}}^{-1})(\bA_{i,j})=\bA_i$. Clearly $\phi$ defines an
isomorphism between $\cA$ and $\cA'$. \qed\epf

\bde\label{dfn:a pointing of a graph of groups}
Let $(\cC,\orG)$ be a graph of groups, a {\em pointing} is a
pair $((\cC',\orG),\delta)$, where $\delta=\{\delta_{i,j}\mid (i,j)\in \orE\}$ is a collection  of elements $\delta_{i,j}\in A_i$ and  $(\cC',\orG)$ is a graph of groups obtained from $(\cC,\orG)$ by
   setting $\alpha'_{i,j}=\ad(\delta_{i,j}^{-1})\after\alpha_{i,j}$, for each $(i,j)\in \orE$.
\ede

\ble\label{lem:all pointings are created equally} Let
$((\cC',\orG),\delta)$ be a pointing of $(\cC,\orG)$. Then
$(\cC',\orG)$ and $(\cC,\orG)$ are isomorphic as graphs of groups.
\ele \bpf In Definition~\ref{dfn:inner morphism} let all $\phi_i$
and $\phi_{i,j}$ be identity maps and let $\delta_{i,j}$ be as
defined in Definition~\ref{dfn:a pointing of a graph of groups}.
This defines an isomorphism $(\cC',\orG)\to (\cC,\orG)$ of graphs of
groups. \qed\epf

\bde\label{dfn:isomorphism of pointings}
An isomorphism between pointings $((\cC,\orG),\delta^{(k)})$ of the concrete graph of groups $(\cC_0,\orG)$ is an inner isomorphism $(\phi,\gamma)$ such that $\gamma=\id$ and there exist $a_i\in \bA_i$ and $a_{i,j}=a_{j,i}\in \bA_{i,j}$ such that $\phi_i=\ad(a_i)$ for each $i\in I$ and $\phi_{i,j}=\ad(a_{i,j})$ for each $(i,j)\in \orE$.
We also require that for each $(i,j)\in \orE$ we have
$ \delta_{i,j}^{(1)}\alpha_{i,j}(a_{i,j})= a_i\delta_{i,j}^{(2)}$.
We then say that the collection $\{a_{i,j}, a_i\mid i\in I, (i,j)\in \orE\}$ induces the isomorphism.
For $\{i,j\}\not\in E$ we then set $a_{i,j}=a_i\ast a_j$.
\ede

\bth\label{thm:CT pointing equivalence} For any admissible Dynkin
diagram $\orG$, there is a natural bijection between the set of
isomorphism classes of concrete CT-structures over $\orG$  and the
set of isomorphism classes of pointings of the concrete graph of
groups $(\cC_0,\orG)$. \eth \bpf Let $\cA$ be a concrete CT
structure over $\orG$.  Then $\cA$ defines a pointing of $\cC_0$ by
setting $\delta_{i,j}=\varphi_{i,j}^{-1}\after\psi_{i,j}$, for each
$(i,j)\in \orE$.

Conversely, given a pointing $((\cC_0,\orG),\delta)$ we define a CT-structure $\cA$ over $\orG$ setting $\varphi_{i,j}=\psi_{i,j}\after\delta_{i,j}^{-1}$ for each $(i,j)\in \orE$.  Of course if $(i,j)\not\in \orE$ then $G_{i,j}=G_i\ast G_j$ so the maps $\varphi_{i,j}$ are the natural ones. The fact that the collection $\{\varphi_{i,j}\mid i,j\in I\}$ defines a concrete CT-structure is immediate since $\ad(\varphi_{i,j})=\ad(\delta_{i,j}^{-1})\after\ad(\psi_{i,j})$, where $\ad(\psi_{i,j})$ takes $\bA_{i,j}$ to $\bA_i$ and $\ad(\delta_{i,j}^{-1})$ preserves $\bA_i$.

We now show that these maps preserve isomorphism classes. First assume that $\cA$ and $\cA'$ are concrete CT structures defined by the collections $\{\varphi_{i,j}\mid i,j\in I\}$ and $\{\varphi_{i,j}'\mid i,j\in I\}$ and $\phi\colon\cA\to\cA'$ is an isomorphism of concrete CT-structures. Fix some $i,j\in I$.
Now we have $\phi_{i,j}\after \varphi_{i,j}= \varphi'_{i,j}\after \phi_i$.  By Lemma~\ref{lem:diagonal automorphism extension} we can assume that, after possibly composing with a diagonal automorphism, $\ad(\phi_i)(\bA_i)=\bA_i$ and since $\cA$ and $\cA'$ are concrete we then also have $\ad(\phi_{i,j})(\bA_{i,j})=\bA_{i,j}$. However the elements of $\Aut_{G_i}(D_i)$ respectively  $\Aut_{G_{i,j}}(G_i, G_j)$ that preserve these complements are exactly the elements of those complements. This means that $\phi_i\in \bA_i$ and $\phi_{i,j}\in \bA_{i,j}$. The collection $\{\phi_{i,j}, \phi_i\mid i,j\in I, (i,j)\in \orE\}$ induces the desired isomorphism between $((\cC_0,\orG),\delta')$ and $((\cC_0,\orG),\delta)$ in the sense of Definition~\ref{dfn:isomorphism of pointings}. Indeed,
$$\delta'_{i,j}\alpha_{i,j}(\phi_{i,j})=((\varphi'_{i,j})^{-1}\psi_{i,j})(\psi_{i,j}^{-1}\phi_{i,j}\psi_{i,j})=(\varphi'_{i,j})^{-1}\phi_{i,j}\psi_{i,j} = \phi_i\varphi_{i,j}^{-1}\psi_{i,j}=\phi_i \delta_{i,j}.$$
Conversely suppose  $\{\phi_{i,j}, \phi_i\mid i,j\in I, (i,j)\in \orE\}$ induces  an isomorphism of pointings $((\cC_0,\orG),\delta')$ and $((\cC_0,\orG),\delta)$. We show that $\phi$ uniquely defines an isomorphism  of CT structures.  Indeed, whenever $(i,j)\in \orE$ we have
$$ \phi_{i,j}\varphi_{i,j}= \phi_{i,j}\psi_{i,j}\delta_{i,j}^{-1}=\psi_{i,j}\alpha_{i,j}(\phi_{i,j})\delta_{i,j}^{-1}=\psi_{i,j}(\delta_{i,j}')^{-1}\phi_i=\varphi_{ij}'\phi_i.$$
In case $(i,j)\not\in \orE$ we simply let
$\phi_{i,j}=\phi_i\ast\phi_j$. \qed\epf

\bco\label{cor:CT pointing equivalence} For any admissible Dynkin
diagram $\orG$, there is a natural bijection between the set of
isomorphism classes of CT-structures over $\orG$  and the set of
isomorphism classes of pointings of the concrete graph of groups
$(\cC_0,\orG)$. \eco \bpf This follows from Corollary~\ref{co:all
amalgams are concrete} and Theorem~\ref{thm:CT pointing
equivalence}. \qed\epf

\subsection{The fundamental group}

\bde\label{def:path group}
For a given graph of groups $(\cC, \orG)$ we define its {\em path group} as follows
$$ \pi(\cC)=((*_{i\in I} A_i )* F(\orE))/ R$$
where $F(\orE)$ is  the free group on the set $\orE$, $*$ denotes free product and $R$ is the following set of relations:  for any $e=(i,j) \in \orE$, we have
\begin{equation}\label{eqn:path group relations}
 \begin{array}{rll}
 e\bar e&=\id &\mbox{ and } \\
 e\cdot \alpha_{\bar{e}}(a)\cdot\bar e &=\alpha_e(a)&\mbox{ for any }a \in  A_e.\\
\end{array}\end{equation}
\ede

\bde\label{def:fundamental group}
Given a graph of groups $(\cC, \orG)$, a {\em path of length $n$} in $\cC$ is a sequence $\gamma=(a_1, e_1, a_2, \ldots, e_{n-1}, a_n)$, where $e_1, \ldots, e_{n-1}$ is an edge path in $\orG$ with vertex sequence $i_1, \ldots, i_n$ and $a_k \in A_{i_k}$ for each $k=1,\ldots,n$. We call $\gamma$ {\em reduced} if it has no returns (i.e.\ $e_{i+1}\ne \bar e_i$ for any $i=1,\ldots,n-2$). The path $\gamma$ defines an element   $|\gamma|=a_1\cdot e_1 \cdot a_2 \cdots e_{n-1} \cdot a_n \in\pi(\cC)$. We denote by $\pi[i,j]$ the collection of elements $|\gamma|$, where $\gamma$ runs through all paths from $i$ to $j$ in $\cC$. Concatenation induces a group operation on $\pi(\cC, i)=\pi[i,i]$ and we call this group {\em the fundamental group of $\cC$ based at $i$}.    \ede

\noindent
From now on the only graph of groups we will consider is the concrete graph $(\cC_0,\orG)$.

\ble \label{lem:move things across a path}Any element  $|\gamma| \in
\pi(\cC_0, i_0)$ can be uniquely realized as $e_1e_2\cdots e_n g$
where  $e_1=(i_0,i_1),\ldots, e_{n-1}=(i_{n-2},i_{n-1})$,
$e_n=(i_{n-1},i_0)$ and $g \in A_{i_0}$. More precisely, if $\gamma
= (e_1, \delta_{i_1},\ldots, \delta_{i_{n-1}}, e_n, \delta_{i_0})$
with $\delta_k\in \bA_k$ for $k=i_0,\ldots,i_{n-1}$, then we have
$g= \beta_{i_1,i_0}(\delta_{i_1})
\beta_{i_2,i_0}(\delta_{i_2})\cdots
\beta_{i_n,i_0}(\delta_{i_n})\delta_{i_0}$. \ele \bpf The first part
is a special case of Corollary 1.13 in~\cite{Bas93} since all maps
$\alpha_{i,j}$ are surjective. The second part follows from the
relations in Definition~\ref{def:path group} and the definition of
$\beta_{l,m}$ preceding Lemma~\ref{lem:walking on trivial paths}.
\qed\epf

\bco\label{compute the fundamental group of c0} $\pi(\cC_0, i_0)\cong
\bA_{i_0}\times \pi(\orG,i_0)$. \eco \bpf By Lemma~\ref{lem:move
things across a path} $\pi(\cC_0, i_0)\cong \bA_{i_0}
\pi(\orG,i_0)$. Also, if $a\in \bA_{i_0}$ and $\gamma
\in\pi(\orG,i_0)$ then $a\cdot\gamma=\gamma \cdot
\beta_{i_0,i_0}(a)=\gamma\cdot a$ so $[\bA_{i_0},\pi(\orG,i_0)]=1$.
Clearly $\bA_{i_0}\cap \pi(\orG,i_0)=1$. \qed\epf

We also need a slight modification of (Corollary 1.10 of
\cite{Bas93}). We first prove the following special case, which uses
the relation~(\ref{eqn:path group relations}).
\ble\label{lem:cor1.10 for returns} If $e\in \orE$ and $\eta=(g_1,
e, g_2, \bar e,  g_3)$ and $\eta'=(g'_1, e, g'_2, \bar e, g'_3)$ are
two paths satisfying  $g_1\alpha_e(\alpha_{\bar
e}^{-1}(g_2))g_3=g'_1\alpha_e(\alpha_{\bar e}^{-1}(g'_2))g'_3$ (so
in particular, $|\eta|=|\eta'|$) then there exist $h_1, h_2\in A_e$
so that $g'_1=g_1\alpha_e(h_1^{-1}), g'_2=\alpha_{\bar
e}(h_1)g_2\alpha_{\bar e}(h_2^{-1}), g_3'= \alpha_e (h_2)g_3$. \ele
\bpf We define $h_1=\alpha_e^{-1}((g_1')^{-1}g_1)$ and
$h_2=\alpha_e^{-1}(g_3'g_3^{-1})$.  The condition on the $g_i$'s can
be rewritten as $\alpha_e(\alpha_{\bar
e}^{-1}(g_2'))=(g_1')^{-1}g_1\alpha_e(\alpha_{\bar
e}^{-1}(g_2))g_3(g_3')^{-1}$. If we apply  $\alpha_e^{-1}$ to this
relation we get $\alpha_{\bar e}^{-1}(g_2')=h_1\alpha_{\bar
e}^{-1}(g_2)h_2^{-1}$. Another application of $\alpha_{\bar e}$
finishes the proof. \qed\epf

We are now ready to prove the following generalization of Corollary 1.10 of \cite{Bas93}.

\bpr\label{prop:homotopy of paths}
Let $\gamma=(g_0,e_1,g_1,\ldots ,e_n,g_n) $ and $\gamma=(g'_0,e_1,g_1',\ldots, e_n,g'_n)$ be two paths with
 $|\gamma|=|\gamma'| $ in $\pi(\cC_0)$. Then there exist elements  $h_i\in A_{e_i}$ ($i=1,2,\ldots,n$) such that
\begin{equation}\label{eqn:prop1}
\begin{array}{rl}
g_0'&=g_0\alpha_{e_1}(h_1^{-1}), \\
g_i' &= \alpha_{\bar e_i}(h_i)g_i\alpha_{e_{i+1}}(h_{i+1}^{-1}),\mbox{ for all }i=1,2,\ldots,n-1,\mbox{ and } \\
g_n'&=\alpha_{\bar e_n}(h_n)g_n.\\
\end{array}
\end{equation}
\epr \bpf If $\gamma$ and $\gamma'$ are reduced then this is just
Corollary 1.10 of \cite{Bas93}. We prove the general case by
induction on the number of returns. Suppose we have a return $e_j
=\bar e_{j+1}$ for some $j=1,\ldots,n$. By "omitting" the return, we
get paths $\dot\gamma=(g_0,e_1,\ldots,e_{j-1},\dot g,
e_{j+2},\ldots,e_n,g_n)$
 and
$\dot\gamma'=(g'_0, e_1,\ldots, e_{j-1},\dot g', e_{j+2},\ldots,
e_n,g'_n)$, where $\dot
g=g_j\alpha_{e_j}(\alpha_{e_{j+1}}^{-1}(g_{j+1}))g_{j+2}$ and $\dot
g'=g_j'\alpha_{e_j}(\alpha_{e_{j+1}}^{-1}(g_{j+1}'))g_{j+2}'$. Using
the relations~(\ref{eqn:path group relations}) we can immediately
see that $|\dot\gamma|=|\gamma|=|\gamma'|=|\dot\gamma'|$. By
induction there exist $h_1,\ldots, h_{j-1}, h_{j+2}, \ldots h_n$
that satisfy the relations~(\ref{eqn:prop1}) for $i\ne  j, j+1, j+2$
as well as the relation $g'= \alpha_{\bar
e_{j-1}}(h_{j-1})g\alpha_{e_{j+2}}(h_{j+2}^{-1})$. We now take $\dot
g_1=  \alpha_{\bar e_{j-1}}(h_{j-1})g_j$, $\dot g'_1=g'_j$, $\dot
g_2 =g_{j+1}$, $\dot g'_2 =g'_{j+1}$, $\dot g'_3= g'_{j+2}$ $\dot
g_3= g_{j+2}\alpha_{e_{j+2}}(h_{j+2}^{-1})$. The paths $(\dot g_1,
e_j,\dot g_2, e_{j+1}, \dot g_3)$ and  $(\dot g'_1, e_j, \dot g'_2,
e_{j+1}, \dot g'_3)$ satisfy  the conditions of
Lemma~\ref{lem:cor1.10 for returns} so there exist $\dot h_1$ and
$\dot h_2$ as in the conclusion of that lemma. Picking $h_j=\dot
h_1$ and $h_{j+1}=\dot h_2$  finishes the proof. \qed\epf

\bde\label{def: fund group of a pointing}
If $((\cC_0, \orG), \delta)$ is a pointing of the graph of groups $(\cC_0, \orG)$ then any  path  $\gamma=e_1 \cdots e_n$ in $\orG$ gives rise to a  path in $\cC_0$ via $\gamma \mapsto \gamma_{\delta}=\delta_{e_1} e_1 \delta_{\bar e_1}^{-1}\delta_{e_2} \cdots e_{n-1}\delta_{\bar e_{n-1}}^{-1}\delta_{e_n}e_n \delta_{\bar e_n}^{-1}$. The map $\gamma\mapsto |\gamma_\delta|$ restricts to a monomorphism  $i_\delta\colon\pi(\orG,i_0) \to \pi(\cC_0,i_0)$. The image of this map is called {\em the fundamental group of the pointing} and denoted by $\pi(\cC_0,i_0,\delta)$. \ede

\ble\label{lem:homomorphism from pi to A} If $((\cC_0, \orG),
\delta)$ is a pointing of $\cC_0$ then there exists a homomorphism
$\Phi\colon\pi(\orG, i_0)\to \bA_{i_0}$ so that
$\pi(\cC_0,i_0,\delta)=\{\gamma\cdot\Phi(\gamma) \mid \gamma \in
\pi(\orG, i_0)\}$. \ele \bpf In view of Corollary~\ref{compute the
fundamental group of c0}, there is a projection homomorphism
$p_{i_0}\colon \pi(\cC_0,i_0)\to \bA_{i_0}$. Now $\Phi=p_{i_0}\after
i_\delta$. The description of the elements in
$\pi(\cC_0,i_0,\delta)$ follows from Lemma~\ref{lem:move things
across a path}. \qed\epf

Note that any two pointings have isomorphic fundamental groups. Therefore the real invariant of the pointing is the actual image of $\pi( \orG, a) \to \pi(\cC_0,a)$ and not its isomorphism class. We in fact have the following.

\bth \label{thm:pointings are fundamental groups} Two pointings of
$\cC_0$ are isomorphic if and only if they have the same fundamental
group. \eth \bpf Suppose that the collection $\{a_{i,j}, a_i\mid i
\in I, (i,j)\in \orE\}$ induces an isomorphism of $((\cC_0, \orG),
\delta)$ to $((\cC_0, \orG), \delta')$. This means that
$a_i^{-1}\delta_{i,j}\alpha_{i,j}(a_{i,j})= \delta_{i,j}'$, for any
$(i,j)\in \orE$. Suppose that $\gamma=e_1, e_2, \ldots, e_n$ is a
path in $\orG$ where without loss of generality we can assume that
$e_i=(i-1,i)$, for each $i=1,\ldots,n$. Then
$$|\gamma_{\delta'}| = \delta'_{0,1} e_1(\delta'_{1,0})^{-1} \delta'_{1,2},e_2, \ldots, e_{n-1}(\delta'_{n-1,n-2})^{-1}\delta'_{n-1,n}e_n (\delta'_{n,n-1})^{-1}$$
Recall the following conditions, for each $(i,j)\in \orE$:
\beq\label{eqn:pointing iso}
a_i^{-1}\delta_{i,j}\alpha_{i,j}(a_{i,j})=\delta_{i,j}' \mbox{ and }
a_{i,j}=a_{j,i}.
\eeq
Since these conditions are met by the collection $\{a_i,a_{i,j}\mid i,j\in I, (i,j)\in \orE\}$ we have
$$\delta'_{i-1,i} e_i(\delta'_{i,i-1})^{-1} =  a_{i-1}^{-1}\delta_{i-1,i}\alpha_{i-1,i}(a_{i-1,i})e_i \alpha_{i,i-1}(a_{i-1,i}^{-1})\delta_{i,i-1}^{-1}a_{i}=a_{i-1}^{-1}\delta_{i-1,i}e_i\delta_{i,i-1}^{-1}a_{i}.$$
Hence if $\gamma$ is a cycle based at $i_0$, then all $a_i$'s cancel and $|\gamma_\delta|=|\gamma_{\delta'}|$ in $\pi(\cC_0,i_0)$. It follows that $\pi(\cC_0,i_0,\delta)=\pi(\cC_0,i_0,\delta')$.

Conversely suppose that $((\cC_0, \orG), \delta)$ and $((\cC_0, \orG), \delta')$ have the same fundamental group.
Recall that we must find a collection $\{a_i,a_{i,j}\mid i,j\in I, (i,j)\in \orE\}$ satisfying the conditions~(\ref{eqn:pointing iso}), for each $(i,j)\in \orE$.
Now to compute $a_j$ for some $j\in I$ we consider a path
 $\gamma=e_1,e_2,\ldots,e_n$ where $e_k=(i_{k-1},i_k)$ for
  $k=1,\ldots,n$ and $i_n=j$.
Now a choice of $a_{i_0}$ uniquely determines  the values of $a_{i_{k}}$ and $a_{i_{k-1}, i_{k}}$, for all $k=1,\ldots,n$, via the conditions~(\ref{eqn:pointing iso}).
Without loss of generality we fix $a_{i_0}=1$.
We claim that the value of $a_j$ computed in this way does not depend on the choice of $\gamma$.

To prove the claim we need to show that if $\gamma$ is a closed path and we set $j=i_0$, then we'll find $a_j=1$ as well.
Note that $\gamma$ is not necessarily simple so that some vertices and edges might be repeated. We ignore this and compute the $a_{i_{k}}$ and $a_{i_{k-1},i_{k}}$ as if they are all distinct. Since $\pi(\cC_0,i_0,\delta)=\pi(\cC_0,i_0,\delta')$ it follows from Lemma~\ref{lem:move things across a path} that $|\gamma_{\delta}|=|\gamma_{\delta'}|$.
Therefore, the paths $\gamma_\delta$ and $\gamma_{\delta'}$ satisfy the conditions of Proposition~\ref{prop:homotopy of paths} with the elements
$g_0=\delta_{e_1}, g_k=\delta_{\bar{e_{k}}}^{-1}\delta_{e_{k+1}}, g_n=\delta_{\bar{e_n}}^{-1}$ respectively $ g'_0=\delta'_{e_1}, g_k'=(\delta'_{\bar{e_{k}}})^{-1}\delta'_{e_{k+1}}, g'_n=(\delta'_{\bar{e_n}})^{-1}$.
On the other hand, by choice of the $a_{e_k}=a_{i_{k-1}, i_k}$ and $a_{i_k}$ we have
$$ \begin{array}{rl}g'_0= &g_0\alpha_{e_1}(a_{e_1}); \\ &\\
g'_k=&\alpha_{\bar {e_k}}(a_{\bar {e_k}}^{-1})\delta_{\bar{e_{k}}}^{-1}a_{i_k}a_{i_k}^{-1}\delta_{e_{k+1}}\alpha_{e_{k+1}}(a_{e_{k+1}})\\ = & \alpha_{\bar {e_k}}(a_{\bar {e_k}}^{-1})\delta_{\bar{e_{k}}}^{-1}\delta_{e_{k+1}}\alpha_{e_{k+1}}(a_{e_{k+1}})\\ = &\alpha_{\bar {e_k}}(a_{\bar {e_k}}^{-1})g_k\alpha_{e_{k+1}}(a_{e_{k+1}});\\ &\\
g'_n=& \alpha_{\bar{e_n}}(a_{\bar{e_n}}^{-1})g_n.
\end{array}$$
In other words, if they exist, the elements $h_k=a_{e_k}^{-1}$
satisfy the conclusion of Proposition~\ref{prop:homotopy of paths}.
In view of the uniqueness of the elements $a_{e_k}$, their
existence thus follows from that proposition.\qed\epf

Theorem~\ref{mainthm:CT structures} is now a consequence of Theorem~\ref{thm:CT pointing equivalence}, Theorem~\ref{thm:pointings are fundamental groups} and Lemma~\ref{lem:homomorphism from pi to A}.

\section{The Curtis-Tits theorem}\label{section:Curtis-Tits Theorem}
Let $G$ be a  simply connected Kac-Moody group that is locally split over a field $\fk$ (in the sense of~\cite{Mu1999}) with an admissible simply laced Dynkin diagram $\Gamma$ over some finite index set $I$.
We shall prove that the Curtis-Tits amalgam for this group is in fact a Curtis-Tits structure for $G$.

Let $(W,\{r_i\}_{i\in I})$ be the Coxeter system of type $\Gamma$.
Then $G$ has a twin BN-pair $(B^+,N,B^-)$ of type $\Gamma$, which gives rise to a Moufang twin-building
 $\De=(\De_+,\De_-,\delta_+,\delta_-,\delta_*)$ of type $\Gamma$, where, for $\vep=\pm$ we have $\De_\vep=G/B^\vep$ and
  $$\begin{array}{lll}
 \delta_\vep (gB^\vep,hB^\vep)&= w \in W&\mbox{ whenever }B^\vep g^{-1}hB^\vep=B^\vep w B^\vep,\\
 \delta_*(gB^+,hB^-)&= w\in W &\mbox{ whenever }B^+g^{-1}hB^+=B^+wB^+.\\
  \end{array}$$
Two chambers $c$ and $d$ are called {\em opposite} if $\delta_*(c,d)=1$.
Fix two opposite chambers $c_+=B^+$ and $c_-=B^-$.
The standard parabolic subgroups of type $(J,\vep)$, where $J\sbe I$ and $\vep=+,-$,
 are the groups $P^\vep_J=B^\vep W_J B^\vep$, where
  $W_J=\langle r_j\mid j\in J\rangle_W$.
The Levi-decomposition of $P_J$ is
 $P_J=U_J\rtimes L_J$, where $L_J$ is called the Levi-component and
  $U_J$ is the unipotent radical of $P_J$.
We shall write $L_i=L_{\{i\}}$ and $L_{i,j}=L_{\{i,j\}}$ for $i,j\in I$.

Since $\Gamma$ is simply laced, condition (co) is satisfied so by~\cite{MuRo1995a}  the local structure of $\De$ determines the global structure. The Curtis-Tits theorem~\cite[Ch. 13]{Tits1974} and \cite{AbrMuh97} yields $G$
 as the universal completion of the following amalgam:
$$\cA=\{L_i,L_{\{i,j\}}\mid i,j\in I\}$$
The fact that $G$ is locally split means that whenever $i$ and $j$ are adjacent, then the $\{i,j\}$-residue on $c_+$ is isomorphic to the building
 associated to the group $\SL_3(\fk)$.
This implies that $L_{i,j}$ is isomorphic to a quotient of $\SL_3(\fk)$ and has
 $\PSL_3(\fk)$ as a quotient.
We call $G$ simply connected if in fact $L_{i,j}\cong \SL_3(\fk)$.
In particular, this means that $L_i\cong \SL_2(\fk)$ and that $L_i$ and $L_j$
 form a standard pair.
Also, whenever $i$ and $j$ are not adjacent in $\Gamma$, $L_i$ and $L_j$ commute so that $L_{i,j}\cong L_i\ast L_j$.
Thus, $\cA$ is a Curtis-Tits structure over $\Gamma$.
The fact that $\cA$ is oriented follows from the observation that for each $i$, the root group $X_i$ of the fundamental positive root ${\alpha_i}$ belongs to $L_i$ and
 $X_i\sbe B^+$.
In particular, $X_i$ and $X_j$ belong to a common Borel group of
 $L_{i,j}$. Thus $\cA$ is an oriented Curtis-Tits structure for $G$.

In the remainder of this section, we shall prove that every oriented CT-structure with admissible Dynkin diagram can be obtained as the Curtis-Tits amalgam of some simply connected Kac-Moody group that is locally split over $\fk$.

Our strategy is as follows. Let $\Gamma$ be an admissible Dynkin diagram and $\cA(\Gamma)$ an oriented CT-structure over some field $\fk$.
The fact that $\cA(\Gamma)$ is oriented allows us to define a Moufang foundation, which by a result of M\"uhlherr is integrable to a twin-building $\Delta$. We then show that if $G$ is the automorphism group of $\Delta$, then the Curtis-Tits amalgam for $G$ is isomorphic to $\cA(\Gamma)$.

\subsection{Sound Moufang foundations and orientable CT-amalgams}\label{sec:Sound Moufang foundations and orientable CT-amalgams}
We shall make use of the following definition of a foundation~\cite{Mu1999}, which is equivalent to the definition in \cite{Tit87}:

\bde Let $\Gamma$ be an admissible Dynkin diagram over $I$.
A {\dfn foundation of type $\Gamma$} is a triple
$$(\{\Delta_{i,j}\mid \{i,j\}\in E\},\{C_{i,j}\mid \{i,j\}\in E\},\{\theta_{j,i,k}\mid \{i,j\},\{i,k\}\in  E\}), $$
satisfying the following conditions:
\begin{itemize}
\item[(Fo1)] $\Delta_{i,j}$ is a building of type $A_2$ for each $\{i,j\}\in E$;
\item[(Fo2)] $C_{i,j}$ is a chamber of $\Delta_{i,j}$ for each $\{i,j\}\in E$;
\item[(Fo3)] $\theta_{j,i,k}$ is a bijection between the $i$-panel on $C_{i,j}$ in $\Delta_{i,j}$
 and the $i$-panel on $C_{i,k}$ in $\Delta_{i,k}$ such that $\theta_{j,i,k}(C_{i,j})=C_{i,k}$ and if
  $i,j,k,l\in I$ are such that $\{i,j\},\{i,k\},\{i,l\}\in E$, then $\theta_{k,i,l}\after\theta_{j,i,k}=\theta_{j,i,l}$.
\end{itemize}
This foundation is said to be of {\em Moufang} type if  $\Delta_{i,j}$ is a Moufang building for each
 $\{i,j\}\in E$ and if  in [Fo3] the map $\theta_{j,i,k}$ induces an isomorphism between the Moufang set
  induced by $\Delta_{i,j}$ on the $i$-panel of $C_{i,j}$ and the Moufang set induced by $\Delta_{i,k}$ on the $i$-panel of $C_{i,k}$.
\ede

We shall now describe how to obtain a Moufang foundation from a given orientable CT-structure $(\cA,\orG)$.
Let $\{X_i,\mid i\in I\}$ be the collection of root groups as in Definition~\ref{dfn:OCT} and let
 $\{B_{i,j}\mid (i,j)\in \orE\}$ be the collection of Borel groups in $G_{i,j}$ such that
  $\varphi_{i,j}(X_i)$ and $\varphi_{j,i}(X_j)$ are contained in $B_{i,j}$ for any $(i,j)\in \orE$ (note that this in fact determines $B_{i,j}$ uniquely).
For each $(i,j)\in \orE$, let $\Delta_{i,j}$ be the Moufang building of type $A_2$ obtained from $G_{i,j}$ via the BN-pair $(B_{i,j},N_{G_{i,j}}(D_{i,j}))$ and let $C_{i,j}$ be the chamber given by $B_{i,j}$.
Now let $i,j,k\in I$ be such that $(i,j),(i,k)\in \orE$.
Let $\cE_{i,j}$ be the element of $G_{i,j}$ given by
$$\left(\begin{array}{@{}ccc@{}}
0 & -1 & 0 \\
1 & 0 & 0 \\ 0 & 0 & 1\\
\end{array}\right)$$ with respect to the ordered basis $\sE_{i,j}$.
Note that the $i$-panel of $\Delta_{i,j}$ containing $C_{i,j}$ is represented by $C_{i,j}$ itself together with the cosets $\varphi_{i,j}(\lambda)\cE_{i,j}B_{i,j}$, where  $\lambda\in {X_{i}}$.
We now define $\theta_{j,i,k}(\varphi_{i,j}(\lambda)\cE_{i,j} B_{i,j})=\varphi_{i,k}(\lambda) \cE_{i,k}B_{i,k}$.
Note that since the structure of the $i$-panel of $\Delta_{i,j}$ on $C_{i,j}$ (resp. of $\Delta_{i,k}$ on $C_{i,k}$) as a Moufang set is entirely determined by $G_i$ the group isomorphism $\theta_{j,i,k}$ preserves this structure. This proves the following.
\ble\label{lem:Moufang foundation from CT structure}
The triple
$$\bF=(\{\Delta_{i,j}\mid \{i,j\}\in E\},\{C_{i,j}\mid \{i,j\}\in E\},\{\theta_{j,i,k}\mid \{i,j\},\{i,k\}\in  E\})$$
obtained from the CT-structure $(\cA,\orG)$ as above, is a Moufang foundation.
\ele

\ble\label{thm:OCT's come from buildings} A CT-structure over an
admissible Dynkin diagram $\Gamma$ is the amalgam coming from the
Curtis-Tits theorem for a twin-building $\Delta$ if and only if it
is orientable. \ele \bpf As proved in the beginning of
Section~\ref{section:Curtis-Tits Theorem}, the amalgam $\cA(\Delta)$
produced by applying the Curtis-Tits theorem to the universal
Kac-Moody group that is an automorphism group for $\Delta$, is an
orientable CT-structure with diagram $\orG$.

Conversely, let $\cA$ be a concrete OCT structure with diagram $\orG$ and $\bF$ be the Moufang foundation constructed from $\cA$ as in Lemma~\ref{lem:Moufang foundation from CT structure}.

Now $\bF$ gives rise to a system $\cK=\{\fk_{i,j},\phi_{i,j}\mid \{i,j\}\in E\}$ as in~\cite[\S 6.5]{Ti1992} in the following way.
Let $(i,j)\in \orE$. By the discussion in loc. cit. we may identify the additive group of $\fk_{i,j}$ with the root group $X_{i}$, which in turn is canonically identified with $\fk$ by viewing $X_{i}$
 as the upper or lower triangular unipotent matrices in $G_i=\SL_2(\fk)$.
The map $\psi_{i,j}$ identifies the field $\fk=\fk_{i,j}$ with the field $\fk$ defining $\Delta_{i,j}$
 and the identification between $\fk_{i,j}$ and $\fk_{j,i}^\circ$ is induced by the
 base change from $\sE_{i,j}$ to $\sE_{j,i}$ which induces the identity on $\fk$.
The inclusion map of the $i$-panel (resp. $j$-panel) on $C_{i,j}$ in $\Delta_{i,j}$ is given by the group isomorphism $\varphi_{ij}$ (resp. $\varphi_{j,i}$) and so the field isomorphism $\phi_{j,i}\colon \fk_{j,i}\to\fk_{i,j}^\circ$ equals $\delta_{i,j}\after\delta_{j,i}^{-1}$. This corresponds to the element
$(\alpha_{j,i}\after\alpha_{i,j}^{-1})(\delta_{i,j})\delta_{j,i}^{-1}$ in $\Aut(\fk_{j,i})$.
Pick a base point $i_0$, for each $(i_0,i)\in \orE$ identify $\fk_{i_0.i}=\fk$ and identify $\Aut(\fk)=\bA_{i_0}$.
Then $\phi_{i,j}$ corresponds to an element of $\Aut(\fk)$ via $\beta_{i_0,i}$. This is how the homomorphism
 $\Phi\colon \pi(\Gamma)\to \Aut(\fk)$ is obtained in loc. cit..
From the definition~\ref{def:fundamental group} and
Lemma~\ref{lem:move things across a path} we see that this
homomorphism coincides with the $\Phi$ defined in
Lemma~\ref{lem:homomorphism from pi to A}. By loc. cit. all sound
Moufang foundations are determined by the homomorphism $\Phi$. This
homomorphism is the same as the homomorphism $\Phi$ defined in
Lemma~\ref{lem:homomorphism from pi to A}. Therefore by
Theorems~\ref{thm:CT pointing equivalence}~and~\ref{thm:pointings
are fundamental groups} every foundation comes from an OCT structure
that is unique up to isomorphism. \qed\epf

\begin{remark}
An alternate proof of Lemma~\ref{thm:OCT's come from buildings}
could be obtained by using the notion of apartments in foundations
as in \cite{Mu1999}.
\end{remark}

\section{Twists of split Kac-Moody groups}\label{section:twists}
\noindent We first note that if the maps $\psi_{ij}$ are as in Definition~\ref{dfn:concrete amalgams}, then the amalgam  $\{G_i,G_{i,j}, \psi_{i,j}\mid i,j\in I, (i,j)\in \orE\}$ has as its universal completion the simply connected split Kac-Moody group ${\mathcal G}_\Gamma(\fk)$ with Dynkin diagram $\Gamma$ over $\fk$. This suggests the following definition.

\bde\label{def:twisted kac_moody group}
For any pointing  $((\cC_0, \orE), \delta)$ of the graph of groups $(\cC_0,\orG)$, the {\dfn $\delta$-twist} is the group ${\mathcal G}_\Gamma^\delta(\fk)$ given by the Curtis-Tits presentation corresponding to $\delta$. More precisely, for any $(i,j)\in \orE$ and $i \in I$ we get a copy $G_i=\SL_2(\fk)$ and $G_{i,j}=G_{j,i}=\SL_3(\fk)$ and the relations given by those in $G_i$ and $G_{i,j}$ together with the following:
\begin{itemize}
\AI{i} if $(i,j)\in\orE$ then $\varphi_{i,j}=\psi_{i,j}\after \delta_{i,j}^{-1}\colon G_i\into G_{i,j}$ identifies $G_i$ with a subgroup of $G_{i,j}$;
\AI{ii} if $(i,j),(j,i)\not\in \orE$, then  $[G_{i}, G_{j}]=1$.
\end{itemize}
\ede

\noindent Corollary~\ref{mainthm:presentation} follows immediately from Theorem~\ref{mainthm:CT structures}.

\noindent We can now make the description of the $\delta$ twists more precise. To that end, let us fix  a spanning tree $\Lambda$ of $\Gamma$, together with a set of directed edges $H$ that does not intersect $\bar{H}=\{\bar e\mid e\in H\}$ and $\Gamma =\Lambda \cup H\cup \bar{H}$. We will construct an amalgam as follows. For each $e\in H$ we take $\delta_e\in \bA_{i_e}$, where $i_e$ is the starting point of $e$. Now let
$$\varphi_e=\left\{\begin{array}{ll}
\psi_e\after \delta_{e}^{-1}&\mbox{ if } e\in H,\\
\psi_e&\mbox{ else. }\\
\end{array}\right.
$$
The resulting amalgam is denoted by $\cA_\delta$.

\bco\label{cor:classification of KM groups} Let $\Gamma$ be a simply
laced Dynkin diagram with no triangles and $\fk$ a field with at
least 4 elements. Any universal Kac-Moody group with diagram
$\Gamma$ that is locally split over $\fk$ is the universal
completion of a unique $\cA_\delta$. \eco \bpf Since the set $H$
corresponds to a unique set of generators for the fundamental group
of $\Gamma$, there is a natural bijection between sets
 $\{\delta_e\mid e\in H\}$ and homomorphisms $\Phi\colon\pi(\orG,i_0)\to \ZZ_2\times \Aut(\fk)$. The result now follows from Theorem~\ref{mainthm:CT structures}.
\qed\epf

\bibliographystyle{elsarticle-num}

%
%\bibliographystyle{plain}

%\bibliography{resref}
\end{document}